\documentclass[12pt,twoside,american,english,british]{article}
\usepackage[T1]{fontenc}
\usepackage[utf8]{inputenc}
\usepackage[a4paper]{geometry}
\usepackage{textcomp}
\usepackage{amsmath}
\usepackage{amsthm}
\usepackage{amssymb}
\usepackage{setspace}
\usepackage{esint}

\makeatletter
\numberwithin{equation}{section}
\numberwithin{figure}{section}
\theoremstyle{plain}
\newtheorem{thm}{\protect\theoremname}[section]
\theoremstyle{plain}
\newtheorem{cor}[thm]{\protect\corollaryname}
\theoremstyle{definition}
\newtheorem{defn}[thm]{\protect\definitionname}
\theoremstyle{plain}
\newtheorem{lem}[thm]{\protect\lemmaname}
\theoremstyle{plain}
\newtheorem{prop}[thm]{\protect\propositionname}
\newcommand{\lyxaddress}[1]{
	\par {\raggedright #1
	\vspace{1.4em}
	\noindent\par}
}

\usepackage[T1]{fontenc}
\usepackage[utf8]{inputenc}
\usepackage[a4paper]{geometry}
\usepackage{color}
\usepackage{dsfont}
\usepackage{amsmath}
\usepackage{amsthm}
\usepackage{amssymb}
\usepackage{setspace}
\usepackage{esint}

\makeatletter
\numberwithin{equation}{section}
\numberwithin{figure}{section}
\theoremstyle{plain}

\@ifundefined{date}{}{\date{}}
\usepackage{babel}
\usepackage{latexsym}
\usepackage{amsthm}
\usepackage{esint}
\usepackage{eucal}
\usepackage{epstopdf}
\usepackage{graphicx}
\usepackage{bigints}
\theoremstyle{plain}

\newtheoremstyle{boldremark}
    {\dimexpr\topsep/2\relax} 
    {\dimexpr\topsep/2\relax} 
    {}          
    {}          
    {\bfseries} 
    {.}         
    {.5em}      
    {}          

\theoremstyle{boldremark}
\newtheorem{brem} [thm] {Remark} 

\usepackage{fancyhdr}
\usepackage{blindtext}
\usepackage{titlesec}

\titleformat{name=\chapter}[display]
{\normalfont\Huge\itshape}
{\titlerule[1pt]\vspace{-33pt}\filleft%
  \parbox[t]{6em}{%
    \raggedleft%
    \rule{\linewidth}{0.5ex}\newline%
    \chaptertitlename\ \thechapter%
  }%
}
{4pc}
{\normalfont\upshape\bfseries\Huge}
\allowdisplaybreaks

\makeatother

\usepackage{babel}
\addto\captionsbritish{\renewcommand{\definitionname}{Definition}}
\addto\captionsbritish{\renewcommand{\lemmaname}{Lemma}}
\addto\captionsbritish{\renewcommand{\theoremname}{Theorem}}
\addto\captionsenglish{\renewcommand{\definitionname}{Definition}}
\addto\captionsenglish{\renewcommand{\lemmaname}{Lemma}}
\addto\captionsenglish{\renewcommand{\theoremname}{Theorem}}
\providecommand{\definitionname}{Definition}
\providecommand{\lemmaname}{Lemma}
\providecommand{\theoremname}{Theorem}

\usepackage[colorlinks,pdfpagelabels,pdfstartview = FitH,bookmarksopen
= true,bookmarksnumbered = true,linkcolor = blue,plainpages =
false,hypertexnames = false,citecolor = blue,pagebackref=false,urlcolor=blue]{hyperref}

\usepackage{amsmath}  
\usepackage{graphicx} 
\def\Yint#1{\mathchoice
    {\YYint\displaystyle\textstyle{#1}}%
    {\YYint\textstyle\scriptstyle{#1}}%
    {\YYint\scriptstyle\scriptscriptstyle{#1}}%
    {\YYint\scriptscriptstyle\scriptscriptstyle{#1}}%
      \!\iint}
\def\YYint#1#2#3{{\setbox0=\hbox{$#1{#2#3}{\iint}$}
    \vcenter{\hbox{$#2#3$}}\kern-.51\wd0}}
\def\longdash{{-}\mkern-3.5mu{-}} 
\def\tiltlongdash{\rotatebox[origin=c]{15}{$\longdash$}}

\def\tiltfiint{\Yint\tiltlongdash}

\makeatother

\usepackage{babel}
\addto\captionsamerican{\renewcommand{\corollaryname}{Corollary}}
\addto\captionsamerican{\renewcommand{\definitionname}{Definition}}
\addto\captionsamerican{\renewcommand{\lemmaname}{Lemma}}
\addto\captionsamerican{\renewcommand{\propositionname}{Proposition}}
\addto\captionsamerican{\renewcommand{\theoremname}{Theorem}}
\addto\captionsbritish{\renewcommand{\corollaryname}{Corollary}}
\addto\captionsbritish{\renewcommand{\definitionname}{Definition}}
\addto\captionsbritish{\renewcommand{\lemmaname}{Lemma}}
\addto\captionsbritish{\renewcommand{\propositionname}{Proposition}}
\addto\captionsbritish{\renewcommand{\theoremname}{Theorem}}
\addto\captionsenglish{\renewcommand{\corollaryname}{Corollary}}
\addto\captionsenglish{\renewcommand{\definitionname}{Definition}}
\addto\captionsenglish{\renewcommand{\lemmaname}{Lemma}}
\addto\captionsenglish{\renewcommand{\propositionname}{Proposition}}
\addto\captionsenglish{\renewcommand{\theoremname}{Theorem}}
\providecommand{\corollaryname}{Corollary}
\providecommand{\definitionname}{Definition}
\providecommand{\lemmaname}{Lemma}
\providecommand{\propositionname}{Proposition}
\providecommand{\theoremname}{Theorem}

\begin{document}
\title{\textbf{Local boundedness for weak solutions to strongly degenerate
orthotropic parabolic equations}}
\author{Pasquale Ambrosio, Simone Ciani}
\date{}

\maketitle
\noindent \textit{\small{}Dedicated to Antonia Passarelli di Napoli
on the occasion of her 60th birthday, celebrated during the ``Workshop
on Variational Problems and PDEs'' (University of Naples “Parthenope”,
June 19-20, 2025).}{\small\par}
\vspace{5mm}
\begin{abstract}
\begin{singlespace}
\noindent We prove the local boundedness of local weak solutions to
the parabolic equation
\[
\partial_{t}u\,=\,\sum_{i=1}^{n}\partial_{x_{i}}\left[(\vert u_{x_{i}}\vert-\delta_{i})_{+}^{p-1}\frac{u_{x_{i}}}{\vert u_{x_{i}}\vert}\right]\,\,\,\,\,\,\,\,\,\,\mathrm{in}\,\,\,\Omega_{T}=\Omega\times(0,T]\,,
\]
where $\Omega$ is a bounded domain in $\mathbb{R}^{n}$ with $n\geq2$,
$p\geq2$, $\delta_{1},\ldots,\delta_{n}$ are non-negative numbers
and $\left(\,\cdot\,\right)_{+}$ denotes the positive part. The main
novelty here is that the above equation combines an orthotropic structure
with a strongly degenerate behavior. The core result of this paper
thus extends a classical boundedness theorem, originally proved for
the parabolic $p$-Laplacian, to a widely degenerate anisotropic setting.
As a byproduct, we also obtain the local boundedness of local weak
solutions to the isotropic counterpart of the above equation.\vspace{0.2cm}
\end{singlespace}
\end{abstract}

\selectlanguage{english}%
\begin{singlespace}

\section{Introduction}
\end{singlespace}

\selectlanguage{british}%
\begin{singlespace}
\noindent $\hspace*{1em}$Let $\Omega$ be a bounded domain in $\mathbb{R}^{n}$,
$n\geq2$, and let $T\in(0,\infty)$. We are interested in the local
boundedness of local weak solutions to the following parabolic equation
\begin{equation}
\partial_{t}u\,=\,\sum_{i=1}^{n}\partial_{x_{i}}\left[(\vert u_{x_{i}}\vert-\delta_{i})_{+}^{p-1}\frac{u_{x_{i}}}{\vert u_{x_{i}}\vert}\right]\,\,\,\,\,\,\,\,\,\,\mathrm{in}\,\,\,\Omega_{T}=\Omega\times(0,T]\,,\label{eq:equation}
\end{equation}
where $p\geq2$, $\delta_{1},\ldots,\delta_{n}$ are non-negative
numbers and $\left(\,\cdot\,\right)_{+}$ stands for the positive
part. Evolutionary equations of this form have been studied since
the 1960s, with significant contributions from the Soviet school;
see, for instance, the work \cite{Vishik} by Vishik. Equation (\ref{eq:equation})
with all $\delta_{i}$ set to zero is also presented explicitly in
several monographs, including \cite{Lions}, \cite[Example 4.A, Chapter III]{Show}
and \cite[Example 30.8]{Zeid}, among others.\\
$\hspace*{1em}$Let us first observe that (\ref{eq:equation}) looks
quite similar to the parabolic $p$-Laplace equation
\begin{equation}
\partial_{t}u\,=\,\sum_{i=1}^{n}\,(\vert Du\vert^{p-2}\,u_{x_{i}})_{x_{i}}\,\,\,\,\,\,\,\,\,\,\mathrm{in}\,\,\,\Omega_{T}\,.\label{eq:parPLAP}
\end{equation}
However, the main novelty of equation (\ref{eq:equation}) is that
it couples the following two features:
\[
\mathit{orthotropic\,\,structure}\,\,\,\,\,\,\,\,\,\,\,\,\,\,\,\mathrm{and}\,\,\,\,\,\,\,\,\,\,\,\,\,\,\,\mathit{strongly\,\,degenerate\,\,behavior}\,.
\]
Indeed, unlike the parabolic $p$-Laplace equation, for which the
loss of ellipticity of the operator $\mathrm{div}(\vert Du\vert^{p-2}Du)$
is restricted to a single point, equation (\ref{eq:equation}) becomes
degenerate on the larger set 
\[
\bigcup_{i=1}^{n}\,\{\vert u_{x_{i}}\vert\leq\delta_{i}\}\,.
\]

\noindent $\hspace*{1em}$A more recent work in which equation (\ref{eq:equation})
appears with all $\delta_{i}$ equal to zero is \cite{BBLVpar}. There,
the authors establish local $L^{\infty}$ estimates for the spatial
gradient of local weak solutions to (\ref{eq:equation}), but confining
their analysis to the case $p\geq2$ and $\max\,\{\delta_{i}\}=0$.
In this special case, as already noted in \cite{BBLVpar}, the basic
regularity theory equally applies to both (\ref{eq:equation}) and
(\ref{eq:parPLAP}). A classical reference in the field is DiBenedetto's
monograph \cite{DiBe}, which provides boundedness results for the
solution $u$ (see \cite[Chapter V]{DiBe}), Hölder continuity estimates
for $u$ (see \cite[Chapter III]{DiBe}), as well as Harnack inequalities
for non-negative solutions (see \cite[Chapter VI]{DiBe}). From a
technical standpoint, there is no distinction to be made between (\ref{eq:parPLAP})
and (\ref{eq:equation}) with all $\delta_{i}$ set to zero. Consequently,
the results in \cite{BBLVpar} and \cite[Chapter V]{DiBe} imply that,
in the case $p\geq2$ and $\max\,\{\delta_{i}\}=0$, the local weak
solutions of (\ref{eq:equation}) are locally Lipschitz continuous
in the spatial variable, uniformly in time.\\
$\hspace*{1em}$The primary goal of this paper is to prove that the
local weak solutions of (\ref{eq:equation}) are locally bounded even
when $\max\,\{\delta_{i}\}>0$, thus extending DiBenedetto's result
\cite[Chapter V, Theorem 4.1]{DiBe} to our anisotropic and more degenerate
setting. More precisely, our main result reads as follows. For notation
and definitions we refer to Section \ref{sec:prelim}.\medskip{}

\end{singlespace}
\begin{thm}
\label{thm:main}Let $n\geq2$ and $p\geq2$. Moreover, assume that
\[
u\,\in\,C_{loc}^{0}\left(0,T;L_{loc}^{2}(\Omega)\right)\cap L_{loc}^{p}\left(0,T;W_{loc}^{1,p}(\Omega)\right)
\]
is a local weak solution of equation $(\ref{eq:equation})$. Then
$u\in L_{loc}^{\infty}(\Omega_{T})$. More precisely,\foreignlanguage{american}{
for every cylinder $[(x_{0},t_{0})+Q(\theta,\rho)]\subset\Omega_{T}$
and every $\sigma\in(0,1)$, we have that:}\\
\foreignlanguage{american}{}\\
\foreignlanguage{american}{$\mathrm{(}a\mathrm{)}$ if $p>2$, the
estimate
\begin{equation}
\underset{[(x_{0},t_{0})\,+\,Q(\sigma\theta,\sigma\rho)]}{\mathrm{ess}\,\sup}\,\vert u\vert\,\leq\,\max\left\{ \rho,\left(\frac{\rho^{p}}{\theta}\right)^{\frac{1}{p-2}},\frac{C}{(1-\sigma)^{\frac{n+p}{2}}}\,\,\sqrt{\frac{\theta}{\rho^{p}}}\,\left(\tiltfiint_{[(x_{0},t_{0})\,+\,Q(\theta,\rho)]}\vert u\vert^{p}\,dx\,dt\right)^{\frac{1}{2}}\right\} \label{eq:superquad}
\end{equation}
holds true for some positive constant $C$ }depending only on $n$,
$p$ and $\max\,\{\delta_{1},\ldots,\delta_{n}\}$;\foreignlanguage{american}{}\\
\foreignlanguage{american}{}\\
\foreignlanguage{american}{$\mathrm{(}b\mathrm{)}$ if $p=2$, the
estimate
\begin{equation}
\underset{[(x_{0},t_{0})\,+\,Q(\sigma\theta,\sigma\rho)]}{\mathrm{ess}\,\sup}\,\vert u\vert\,\leq\,\max\left\{ \rho,\frac{C}{(1-\sigma)^{\frac{n+2}{2}}}\,\,\sqrt{\left(\frac{\rho^{2}}{\theta}\right)^{\frac{n}{2}}+\,\frac{\theta}{\rho^{2}}}\,\left(\tiltfiint_{[(x_{0},t_{0})\,+\,Q(\theta,\rho)]}\vert u\vert^{2}\,dx\,dt\right)^{\frac{1}{2}}\right\} \label{eq:quad}
\end{equation}
holds true for some positive constant $C$ }depending only on $n$
and $\max\,\{\delta_{1},\ldots,\delta_{n}\}$.
\end{thm}

\noindent $\hspace*{1em}$The proof of Theorem \ref{thm:main} relies
on an adaptation of De Giorgi's iteration technique; see, for instance,
\cite{LadSolUra} for the nondegenerate case. In Section \ref{sec:energy_estimate},
we first establish a local energy estimate (Proposition \ref{prop:PropEnergy}),
which allows us to control the superlevel sets of the local weak solution
$u$ through suitable cut-off functions. We then construct a sequence
of shrinking cylinders $Q_{j}$ and increasing levels $k_{j}>0$,
and use the energy estimate to derive recursive inequalities for the
integral quantities
\[
Y_{j}:=\,\tiltfiint_{Q_{j}}(u-k_{j})_{+}^{p}\,dx\,dt\,.
\]
These inequalities fall within the scope of a general iteration lemma
(Lemma \ref{lem:Giusti}), ensuring that $Y_{j}\to0$ as $j\to\infty$.
Consequently, the measure of the limiting superlevel set of $u$ vanishes,
which yields the local boundedness of $u$ and the local $L^{\infty}$
estimates (\ref{eq:superquad})$-$(\ref{eq:quad}).\\
$\hspace*{1em}$It is worth recalling that, in the elliptic setting,
the local boundedness of solutions to anisotropic problems has been
extensively investigated, as can be seen, for instance, in \cite{AmbCupMas,BiaCupMas2,Cianchi,CuMaMa1,CuMaMa2,CuMaMa3,CuMaMa4,FeoPasPos,FusSbo1,FusSbo2,Korolev,LeoMas,Strof}.
Pioneering contributions are due to Kolod\={\i}\u{\i} \cite{Kolod}.
We further mention the more recent work by DiBenedetto, Gianazza and
Vespri \cite{DiBeGiaVes}, where precise a priori $L^{\infty}$ estimates
for the solutions are established (see Section 6 there). For an insight
into the parabolic anisotropic setting, we recall that the local boundedness
of weak solutions is indeed a direct consequence of Caccioppoli's
estimates; see, for example, \cite[Section 6]{Ciani} and the references
therein.\\
$\hspace*{1em}$In the final part of this paper, we will focus on
the local weak solutions of the parabolic equation
\begin{equation}
\partial_{t}u-\mathrm{div}\left((\vert Du\vert-\lambda)_{+}^{p-1}\frac{Du}{\vert Du\vert}\right)=0\,\,\,\,\,\,\,\,\,\,\mathrm{in}\,\,\,\Omega_{T}\,,\label{eq:AmbPass}
\end{equation}
where $p\geq2$ and $\lambda>0$ is a fixed parameter. Indeed, by
slightly modifying the proof of Theorem \ref{eq:equation}, we can
establish the following $L^{\infty}$-regularity result.
\begin{cor}
\label{cor:corollario}Let $n\geq2$, $p\geq2$ and $\lambda>0$.
Moreover, assume that 
\[
u\,\in\,C_{loc}^{0}\left(0,T;L_{loc}^{2}(\Omega)\right)\cap L_{loc}^{p}\left(0,T;W_{loc}^{1,p}(\Omega)\right)
\]
is a local weak solution of equation $(\ref{eq:AmbPass})$. Then $u\in L_{loc}^{\infty}(\Omega_{T})$.
More precisely,\foreignlanguage{american}{ for every cylinder $[(x_{0},t_{0})+Q(\theta,\rho)]\subset\Omega_{T}$
and every $\sigma\in(0,1)$, we have that:}\\
\foreignlanguage{american}{}\\
\foreignlanguage{american}{$\mathrm{(}a\mathrm{)}$ if $p>2$, the
estimate
\begin{equation}
\underset{[(x_{0},t_{0})\,+\,Q(\sigma\theta,\sigma\rho)]}{\mathrm{ess}\,\sup}\,\vert u\vert\,\leq\,\max\left\{ \rho,\left(\frac{\rho^{p}}{\theta}\right)^{\frac{1}{p-2}},\frac{C}{(1-\sigma)^{\frac{n+p}{2}}}\,\,\sqrt{\frac{\theta}{\rho^{p}}}\,\left(\tiltfiint_{[(x_{0},t_{0})\,+\,Q(\theta,\rho)]}\vert u\vert^{p}\,dx\,dt\right)^{\frac{1}{2}}\right\} \label{eq:superquad_2}
\end{equation}
holds true for some positive constant $C$ }depending only on $n$,
$p$ and $\lambda$;\foreignlanguage{american}{}\\
\foreignlanguage{american}{}\\
\foreignlanguage{american}{$\mathrm{(}b\mathrm{)}$ if $p=2$, the
estimate
\begin{equation}
\underset{[(x_{0},t_{0})\,+\,Q(\sigma\theta,\sigma\rho)]}{\mathrm{ess}\,\sup}\,\vert u\vert\,\leq\,\max\left\{ \rho,\frac{C}{(1-\sigma)^{\frac{n+2}{2}}}\,\,\sqrt{\left(\frac{\rho^{2}}{\theta}\right)^{\frac{n}{2}}+\,\frac{\theta}{\rho^{2}}}\,\left(\tiltfiint_{[(x_{0},t_{0})\,+\,Q(\theta,\rho)]}\vert u\vert^{2}\,dx\,dt\right)^{\frac{1}{2}}\right\} \label{eq:quad_2}
\end{equation}
holds true for some positive constant $C$ }depending only on $n$
and $\lambda$. 
\end{cor}

\noindent $\hspace*{1em}$When $\lambda>0$, the main feature of equation
(\ref{eq:AmbPass}) is that it exhibits a strong degeneracy, coming
from the fact that its modulus of ellipticity vanishes in the region
$\{\vert Du\vert\leq\lambda\}$, and hence its principal part behaves
like a $p$-Laplace operator only for large values of $\vert Du\vert$.\\
$\hspace*{1em}$The gradient regularity of weak solutions to (\ref{eq:AmbPass})
has been recently studied in \cite{Amb1,Amb2,AmbBau,AmbPass,BoDuGiPa,GenPas}.
In particular, in \cite{AmbBau} the authors establish local $L^{\infty}$
bounds for the spatial gradient of solutions to equations and systems
of the form (\ref{eq:AmbPass}) in the whole range $p>1$. Therefore,
by combining the results of \cite{AmbBau} with Corollary \ref{cor:corollario},
one obtains that local weak solutions of (\ref{eq:AmbPass}) are locally
Lipschitz continuous in the spatial variable, uniformly in time.

\noindent \begin{brem}If the intrinsic relation $\theta=\rho^{p}$
is imposed, the explicit geometric ratios on the right-hand side of
(\ref{eq:superquad}), (\ref{eq:quad}), (\ref{eq:superquad_2}) and
(\ref{eq:quad_2}) reduce to dimension-free constants and, in this
sense, the aforementioned estimates are ``dimensionless”. However,
if $\theta=\rho^{p}$, (\ref{eq:superquad}) and (\ref{eq:superquad_2})
are not homogeneous in $u$.\end{brem}

\subsection{Plan of the paper }

$\hspace*{1em}$The paper is organized as follows. Section \ref{sec:prelim}
is devoted to the preliminaries: after a list of classical notations
and some essential lemmas, we recall the basic properties of Steklov
averages. In Section \ref{sec:energy_estimate}, we establish a local
energy estimate for the local weak solutions of (\ref{eq:equation}).
This estimate is then used in Section \ref{sec:IterIneq} to derive
a family of local iterative inequalities. In Section \ref{sec:main proofs},
these inequalities are employed to complete the proof of Theorem \ref{eq:equation}.
Finally, in the same section, we include the proof of Corollary \ref{cor:corollario}
for completeness.
\selectlanguage{english}%
\begin{singlespace}

\section{Notation and preliminaries\label{sec:prelim}}
\end{singlespace}

\selectlanguage{british}%
\noindent $\hspace*{1em}$In this paper we shall denote by $C$ or
$c$ a general positive constant that may vary on different occasions,
even within the same line of estimates. Relevant dependencies on parameters
and special constants will be suitably emphasized using parentheses
or subscripts. The norm we use on $\mathbb{R}^{n}$ will be the standard
Euclidean one and it will be denoted by $\left|\,\cdot\,\right|$.
In particular, for the vectors $\xi,\eta\in\mathbb{R}^{n}$, we write
$\langle\xi,\eta\rangle$ for the usual inner product and $\left|\xi\right|:=\langle\xi,\xi\rangle^{\frac{1}{2}}$
for the corresponding Euclidean norm.\\
$\hspace*{1em}$In what follows, we use the notation 
\[
K_{\rho}:=\,(-\rho,\rho)^{n},\,\,\,\,\,\,\,\,\,\,\,\,\rho>0,
\]
for the $n$-dimensional open cube centered at the origin with side
length $2\rho$. If $x_{0}\in\mathbb{R}^{n}$, we denote by $[x_{0}+K_{\rho}]$
the cube of center $x_{0}$ and side length $2\rho$ which is congruent
to $K_{\rho}$, i.e., 
\[
[x_{0}+K_{\rho}]:=\left\{ x\in\mathbb{R}^{n}:\underset{1\,\leq\,i\,\leq\,n}{\max}\,\vert x_{i}-x_{0,i}\vert<\rho\right\} .
\]
Moreover, for a positive number $\theta$, we consider the cylinder
\[
Q(\theta,\rho):=\,K_{\rho}\times(-\theta,0)\,,
\]
and if $(x_{0},t_{0})\in\mathbb{R}^{n+1}$, we let $[(x_{0},t_{0})+Q(\theta,\rho)]$
denote the cylinder with vertex at $(x_{0},t_{0})$ congruent to $Q(\theta,\rho)$,
i.e.,
\[
[(x_{0},t_{0})+Q(\theta,\rho)]:=\,[x_{0}+K_{\rho}]\times(t_{0}-\theta,t_{0})\,.
\]
For a general cylinder $Q=B\times(t_{0},t_{1})$, where $B\subset\mathbb{R}^{n}$
and $t_{0}<t_{1}$, we denote by 
\[
\partial_{\mathrm{par}}Q:=\,(\overline{B}\times\{t_{0}\})\cup(\partial B\times(t_{0},t_{1}))
\]
the usual \textit{parabolic boundary} of $Q$, which is nothing but
its standard topological boundary without the upper cap $\overline{B}\times\{t_{1}\}$.

\selectlanguage{english}%
\noindent $\hspace*{1em}$If $E\subseteq\mathbb{R}^{k}$ is a Lebesgue-measurable
set, then we will denote by $\vert E\vert$ its $k$-dimensional Lebesgue
measure. \foreignlanguage{american}{When $0<\vert E\vert<\infty$,
the mean value of a function $v\in L^{1}(E)$ is defined by 
\[
\fint_{E}v(y)\,dy\,:=\,\frac{1}{\vert E\vert}\int_{E}v(y)\,dy\,.
\]
}$\hspace*{1em}$Now \foreignlanguage{british}{let $F:\mathbb{R}^{n}\to\mathbb{R}$
and $G:\mathbb{R}^{n}\to\mathbb{R}$ be the functions defined respectively
by
\begin{equation}
F(\xi):=\sum_{i=1}^{n}\,\frac{1}{p}\,(\vert\xi_{i}\vert-\delta_{i})_{+}^{p}\,\,\,\,\,\,\,\,\,\,\,\,\,\,\,\,\,\mathrm{and}\,\,\,\,\,\,\,\,\,\,\,\,\,\,\,\,\,G(\xi):=\,\frac{1}{p}\,(\vert\xi\vert-\lambda)_{+}^{p}\,.\label{eq:FG}
\end{equation}
}In this work, we define a local weak solution of (\ref{eq:equation})
and of (\ref{eq:AmbPass}) as follows.\medskip{}

\selectlanguage{british}%
\begin{defn}
\noindent A function $u$ is a \textit{local weak solution} of equation
(\ref{eq:equation}) if 
\[
u\,\in\,C_{loc}^{0}\left(0,T;L_{loc}^{2}(\Omega)\right)\cap L_{loc}^{p}\left(0,T;W_{loc}^{1,p}(\Omega)\right)
\]
and, for every test function $\varphi\in C_{0}^{\infty}(\Omega_{T})$,
\begin{equation}
\iint_{\Omega_{T}}\left(u\,\partial_{t}\varphi\,-\langle D_{\xi}F(Du),D\varphi\rangle\right)dx\,dt\,=\,0\,.\label{eq:locweaksol}
\end{equation}
\end{defn}

\begin{defn}
\noindent A function $u$ is a \textit{local weak solution} of equation
(\ref{eq:AmbPass}) if 
\[
u\,\in\,C_{loc}^{0}\left(0,T;L_{loc}^{2}(\Omega)\right)\cap L_{loc}^{p}\left(0,T;W_{loc}^{1,p}(\Omega)\right)
\]
and, for every test function $\varphi\in C_{0}^{\infty}(\Omega_{T})$,
\begin{equation}
\iint_{\Omega_{T}}\left(u\,\partial_{t}\varphi\,-\langle D_{\xi}G(Du),D\varphi\rangle\right)dx\,dt\,=\,0\,.\label{eq:locweaksol_2}
\end{equation}
\end{defn}

\noindent $\hspace*{1em}$We now recall some tools that will be useful
to prove our results. We begin with the following interpolation inequality,
whose proof can be found in \cite[Proposition 3.1, Chapter I]{DiBe}.
\begin{lem}
\label{lem:interpolation}Let $n\geq2$ and $1\leq r,s<\infty$. Then,
there exists a positive constant $C$, depending only on $n$, $r$
and $s$, such that for every $v\in L^{\infty}\left(0,T;L^{r}(\Omega)\right)\cap L^{s}\left(0,T;W_{0}^{1,s}(\Omega)\right)$
we have 
\[
{\displaystyle \iint_{\Omega_{T}}\vert v(x,t)\vert^{q}\,dx\,dt\,\leq\,C^{q}\left(\iint_{\Omega_{T}}{\displaystyle \vert Dv(x,t)\vert^{s}}\,dx\,dt\right)\left(\underset{0\,<\,t\,<\,T}{\mathrm{ess}\,\sup}\int_{\Omega}\vert v(x,t)\vert^{r}\,dx\right)^{\frac{s}{n}}},
\]
where ${\displaystyle q=s\,\frac{n+r}{n}}$.
\end{lem}

\noindent $\hspace*{1em}$The next lemma is crucial to establish our
main result; see \cite[Lemma 7.1]{Giu} for a proof.
\begin{lem}
\label{lem:Giusti}Let $\alpha>0$ and let $\{Y_{j}\}_{j\,\in\,\mathbb{N}_{0}}$
be a sequence of positive real numbers, satisfying the recursive inequalities
\[
Y_{j+1}\,\leq\,C\,b^{j}\,Y_{j}^{1+\alpha}
\]
where $C>0$ and $b>1$. If $Y_{0}\leq C^{-\,\frac{1}{\alpha}}\,b^{-\,\frac{1}{\alpha^{2}}}$,
then 
\[
\lim_{j\to\infty}Y_{j}=0\,.
\]
\end{lem}

\selectlanguage{american}%

\subsection{Steklov averages}

\selectlanguage{british}%
$\hspace*{1em}$\foreignlanguage{american}{In this section, we recall
the definition and some elementary properties of Steklov averages.
Let us denote a domain in space-time by $Q':=\Omega'\times I$, where
$\Omega'\subseteq\Omega$ is a bounded domain and $I:=(t_{0},t_{1})\subseteq(0,T)$.
For every $h\in(0,t_{1}-t_{0})$ and $v\in L^{1}(\Omega'\times I,\mathbb{R}^{k})$,
where $k\in\mathbb{N}$, the \textit{Steklov average} $[v]_{h}(\cdot,t)$
is defined by 
\[
[v]_{h}(x,t):=\begin{cases}
\begin{array}{cc}
{\displaystyle \frac{1}{h}\int_{t}^{t+h}v(x,s)\,ds} & \,\,\,\text{if }\,t\in(t_{0},t_{1}-h],\\
0 & \,\,\,\text{if }\,t\in(t_{1}-h,t_{1}),
\end{array}\end{cases}
\]
for $x\in\Omega'$. This definition implies, for almost every $(x,t)\in\Omega'\times(t_{0},t_{1}-h)$,
\[
\frac{\partial[v]_{h}}{\partial t}(x,t)\,=\,\frac{v(x,t+h)-v(x,t)}{h}\,.
\]
}

\noindent $\hspace*{1em}$\foreignlanguage{american}{The proof of
the following result is straightforward from the theory of Lebesgue
spaces (see \cite[Lemma 3.2, Chapter I]{DiBe}).\medskip{}
}
\selectlanguage{american}%
\begin{lem}
\label{lem:Stek}Let $q,r\geq1$ and $v\in L^{r}\left(t_{0},t_{1};L^{q}(\Omega')\right)$.
Then, as $h\rightarrow0$, $[v]_{h}$ converges to $v$ in $L^{r}\left(t_{0},t_{1}-\varepsilon;L^{q}(\Omega')\right)$
for every $\varepsilon\in(0,t_{1}-t_{0})$. If $v\in C^{0}\left(t_{0},t_{1};L^{q}(\Omega')\right)$,
then as $h\rightarrow0$, $[v]_{h}(\cdot,t)$ converges to $v(\cdot,t)$
in $L^{q}(\Omega')$ for every $t\in(t_{0},t_{1}-\varepsilon)$, $\forall\,\,\varepsilon\in(0,t_{1}-t_{0})$.
\end{lem}

\selectlanguage{british}%
\noindent $\hspace*{1em}$\foreignlanguage{american}{A very useful
formulation, equivalent to (\ref{eq:locweaksol}), is the one involving
Steklov averages. Assume that $u\in C_{loc}^{0}\left(0,T;L_{loc}^{2}(\Omega)\right)\cap L_{loc}^{p}\left(0,T;W_{loc}^{1,p}(\Omega)\right)$
is a local weak solution of (\ref{eq:equation}) in $\Omega_{T}$
and let $h\in(0,T)$. Then, the Steklov average $[u]_{h}$ satisfies
\begin{equation}
\int_{\mathcal{K}\times\{\tau\}}\left(\frac{\partial[u]_{h}}{\partial t}\cdot\varphi\,+\langle[D_{\xi}F(Du)]_{h},D\varphi\rangle\right)dx\,=\,0\label{eq:Steklov}
\end{equation}
for every compact subset $\mathcal{K}$ of $\Omega$, for all $\tau\in(0,T-h]$
and all test functions 
\[
\varphi\,\in\,C_{loc}^{0}\left(0,T;L^{2}(\mathcal{K})\right)\cap L_{loc}^{p}\left(0,T;W_{0}^{1,p}(\mathcal{K})\right).
\]
If $u$ is a local weak solution of (\ref{eq:AmbPass}), then the
Steklov average formulation of (\ref{eq:locweaksol_2}) is obtained
by simply replacing $D_{\xi}F(Du)$ with $D_{\xi}G(Du)$ in (\ref{eq:Steklov}).}

\section{A local energy estimate\label{sec:energy_estimate} }

\noindent $\hspace*{1em}$The proof of Theorem \ref{thm:main} is
based on the following local energy estimate. Throughout this section
and the sequel, $(x_{0},t_{0})\in\Omega_{T}$ and $\theta,\rho>0$
are such that $[(x_{0},t_{0})+Q(\theta,\rho)]\subset\Omega_{T}$,
while $\zeta$ denotes a piecewise smooth cut-off function in $[(x_{0},t_{0})+Q(\theta,\rho)]$
satisfying
\[
0\leq\zeta\leq1,\,\,\,\,\,\,\,\,\,\,\,\,\,\Vert D\zeta\Vert_{\infty}<+\infty,\,\,\,\,\,\,\,\,\,\,\,\,\,\zeta\equiv0\,\,\,\,\,\,\mathrm{on\,\,}\mathrm{the}\,\,\mathrm{parabolic}\,\,\mathrm{boundary}\,\,\mathrm{of}\,\,[(x_{0},t_{0})+Q(\theta,\rho)]\,.
\]

\begin{prop}
\noindent \label{prop:PropEnergy}Let $n\geq2$ and $p\geq2$. Moreover,
assume that $u$ is a local weak solution of equation $(\ref{eq:equation})$.
Then, for every level $k>0$ we have\begin{align}\label{eq:energy_est}
&\underset{t_{0}-\theta\,<\,\tau\,<\,t_{0}}{\sup}\,\int_{[x_{0}+K_{\rho}]}(u-k)_{+}^{2}\,\zeta^{p}(x,\tau)\,dx\,+\sum_{i=1}^{n}\iint_{[(x_{0},t_{0})+Q(\theta,\rho)]}(\vert u_{x_{i}}\vert-\delta_{i})_{+}^{p}\,\zeta^{p}\,\mathds{1}_{\{u\,>\,k\}}\,dx\,dt\nonumber\\
&\,\,\,\,\,\,\,\leq\,p\iint_{[(x_{0},t_{0})+Q(\theta,\rho)]}(u-k)_{+}^{2}\,\zeta^{p-1}\,\partial_{t}\zeta\,dx\,dt\,+\,C\iint_{[(x_{0},t_{0})+Q(\theta,\rho)]}(u-k)_{+}^{p}\,\vert D\zeta\vert^{p}\,dx\,dt\,
\end{align}for a positive constant $C$ depending only on $n$ and $p$.
\end{prop}

\noindent \begin{proof}[\bfseries{Proof}]After a translation, we
may assume that $(x_{0},t_{0})=(0,0)$. Hence, it suffices to prove
\eqref{eq:energy_est} for the cylinder $Q(\theta,\rho)$. In (\ref{eq:Steklov})
we take the test functions 
\[
\varphi=([u]_{h}-k)_{+}\,\zeta^{p}
\]
and integrate with respect to time over $(-\theta,\tau)$, with $\tau\in(-\theta,0)$.
We thus obtain 
\begin{equation}
\int_{-\theta}^{\tau}\int_{K_{\rho}}\frac{\partial[u]_{h}}{\partial t}\,([u]_{h}-k)_{+}\,\zeta^{p}\,dx\,dt\,+\iint_{Q^{\tau}}\langle[A(Du)]_{h},D[([u]_{h}-k)_{+}\,\zeta^{p}]\rangle\,dx\,dt\,=\,0\,,\label{eq:weak1}
\end{equation}
where, for convenience of notation, we have set 
\[
Q^{\tau}:=\,K_{\rho}\times(-\theta,\tau)\,\,\,\,\,\,\,\,\,\,\,\,\,\,\,\,\mathrm{and}\,\,\,\,\,\,\,\,\,\,\,\,\,\,\,A(\eta):=\,D_{\xi}F(\eta)\,,\,\,\,\,\,\eta\in\mathbb{R}^{n}.
\]
The first term in (\ref{eq:weak1}) can be rewritten as
\[
\int_{-\theta}^{\tau}\int_{K_{\rho}}\frac{\partial[u]_{h}}{\partial t}\,([u]_{h}-k)_{+}\,\zeta^{p}\,dx\,dt\,=\,\frac{1}{2}\int_{-\theta}^{\tau}\int_{K_{\rho}}\frac{\partial([u]_{h}-k)_{+}^{2}}{\partial t}\,\zeta^{p}\,dx\,dt\,.
\]
Therefore, integrating by parts, using that $\zeta\equiv0$ on $\partial_{\mathrm{par}}Q(\theta,\rho)$
and letting $h\to0$, by Lemma \ref{lem:Stek} we have\begin{align}\label{eq:limite1}
&\int_{-\theta}^{\tau}\int_{K_{\rho}}\frac{\partial[u]_{h}}{\partial t}\,([u]_{h}-k)_{+}\,\zeta^{p}\,dx\,dt\,\longrightarrow\nonumber\\
&\,\,\,\,\,\,\,\frac{1}{2}\int_{K_{\rho}}(u-k)_{+}^{2}\,\zeta^{p}(x,\tau)\,dx\,-\,\frac{p}{2}\iint_{Q^{\tau}}(u-k)_{+}^{2}\,\zeta^{p-1}\,\partial_{t}\zeta\,dx\,dt\,.
\end{align} Now observe that $\vert A(Du)\vert\leq\sqrt{n}\,\vert Du\vert^{p-1}$.
Then, taking the limit as $h\to0$ in the second term of (\ref{eq:weak1}),
we can apply Lemma \ref{lem:Stek} again. Thus we get\begin{align}\label{eq:limite2}
&\iint_{Q^{\tau}}\langle[A(Du)]_{h},D[([u]_{h}-k)_{+}\,\zeta^{p}]\rangle\,dx\,dt\,\longrightarrow\nonumber\\
&\,\,\,\,\,\,\,\iint_{Q^{\tau}}\langle A(Du),D(u-k)_{+}\rangle\,\zeta^{p}\,dx\,dt\,+\,p\iint_{Q^{\tau}}\langle A(Du),D\zeta\rangle\,(u-k)_{+}\,\zeta^{p-1}\,dx\,dt\,.
\end{align}From (\ref{eq:weak1})$-$\eqref{eq:limite2}, we then obtain\begin{align*}
&\frac{1}{2}\int_{K_{\rho}}(u-k)_{+}^{2}\,\zeta^{p}(x,\tau)\,dx\,+\iint_{Q^{\tau}}\langle A(Du),D(u-k)_{+}\rangle\,\zeta^{p}\,dx\,dt\\
&\,\,\,\,\,\,\,=\,\frac{p}{2}\iint_{Q^{\tau}}(u-k)_{+}^{2}\,\zeta^{p-1}\,\partial_{t}\zeta\,dx\,dt\,-\,p\iint_{Q^{\tau}}\langle A(Du),D\zeta\rangle\,(u-k)_{+}\,\zeta^{p-1}\,dx\,dt\,.
\end{align*}We now estimate\begin{align*}
\iint_{Q^{\tau}}\langle A(Du),D(u-k)_{+}\rangle\,\zeta^{p}\,dx\,dt\,&=\,\sum_{i=1}^{n}\iint_{Q^{\tau}\,\cap\,\{u\,>\,k\}}(\vert u_{x_{i}}\vert-\delta_{i})_{+}^{p-1}\,\vert u_{x_{i}}\vert\,\zeta^{p}\,dx\,dt\\
&\geq\,\sum_{i=1}^{n}\iint_{Q^{\tau}\,\cap\,\{u\,>\,k\}}(\vert u_{x_{i}}\vert-\delta_{i})_{+}^{p}\,\zeta^{p}\,dx\,dt\,.
\end{align*}Furthermore, applying Young's inequality with $\varepsilon>0$, we
have\begin{align*}
&-\,p\iint_{Q^{\tau}}\langle A(Du),D\zeta\rangle\,(u-k)_{+}\,\zeta^{p-1}\,dx\,dt\\
&\,\,\,\,\,\,\,\leq\,p\sum_{i=1}^{n}\iint_{Q^{\tau}}(\vert u_{x_{i}}\vert-\delta_{i})_{+}^{p-1}\,\vert\zeta_{x_{i}}\vert\,(u-k)_{+}\,\zeta^{p-1}\,dx\,dt\\
&\,\,\,\,\,\,\,\leq\,\varepsilon(p-1)\sum_{i=1}^{n}\iint_{Q^{\tau}\,\cap\,\{u\,>\,k\}}(\vert u_{x_{i}}\vert-\delta_{i})_{+}^{p}\,\zeta^{p}\,dx\,dt\,+\,\frac{1}{\varepsilon^{p-1}}\sum_{i=1}^{n}\iint_{Q^{\tau}}\vert\zeta_{x_{i}}\vert^{p}\,(u-k)_{+}^{p}\,dx\,dt\,.
\end{align*}Choosing $\varepsilon=\frac{1}{2(p-1)}$ and collecting the three
previous estimates, we get\begin{align*}
&\int_{K_{\rho}}(u-k)_{+}^{2}\,\zeta^{p}(x,\tau)\,dx\,+\sum_{i=1}^{n}\iint_{Q^{\tau}\,\cap\,\{u\,>\,k\}}(\vert u_{x_{i}}\vert-\delta_{i})_{+}^{p}\,\zeta^{p}\,dx\,dt\\
&\,\,\,\,\,\,\,\leq\,p\iint_{Q^{\tau}}(u-k)_{+}^{2}\,\zeta^{p-1}\,\partial_{t}\zeta\,dx\,dt\,+C\iint_{Q^{\tau}}(u-k)_{+}^{p}\,\vert D\zeta\vert^{p}\,dx\,dt\,,
\end{align*}where $C$ is a positive constant depending only on $n$ and $p$.
Recalling that $\tau\in(-\theta,0)$ is arbitrary, from the above
inequality we obtain\begin{align*}
&\underset{-\theta\,<\,\tau\,<\,0}{\sup}\,\int_{K_{\rho}}(u-k)_{+}^{2}\,\zeta^{p}(x,\tau)\,dx\,+\sum_{i=1}^{n}\iint_{Q(\theta,\rho)}(\vert u_{x_{i}}\vert-\delta_{i})_{+}^{p}\,\zeta^{p}\,\mathds{1}_{\{u\,>\,k\}}\,dx\,dt\\
&\,\,\,\,\,\,\,\leq\,p\iint_{Q(\theta,\rho)}(u-k)_{+}^{2}\,\zeta^{p-1}\,\partial_{t}\zeta\,dx\,dt\,+\,C\iint_{Q(\theta,\rho)}(u-k)_{+}^{p}\,\vert D\zeta\vert^{p}\,dx\,dt\,.
\end{align*}This concludes the proof.\end{proof}

\section{Local iterative inequalities\label{sec:IterIneq}}

\noindent $\hspace*{1em}$An essential ingredient in the proof of
Theorem \ref{thm:main} is a family of iterative inequalities. We
shall now derive them, starting from the energy estimate \eqref{eq:energy_est}.
After a translation, we may assume that $(x_{0},t_{0})$ coincides
with the origin. Fixed $\sigma\in(0,1)$, we consider the sequences
\[
\rho_{j}:=\,\sigma\rho\,+\,\frac{(1-\sigma)}{2^{j}}\,\rho\,,\,\,\,\,\,\,\,\,\,\,\,\,\,\theta_{j}:=\,\sigma\theta\,+\,\frac{(1-\sigma)}{2^{j}}\,\theta\,,\,\,\,\,\,\,\,\,\,\,j\in\mathbb{N}_{0}\,,
\]
and the corresponding cylinders $Q_{j}:=Q(\theta_{j},\rho_{j})$.
From the definitions it follows that 
\[
Q_{0}\,=\,Q(\theta,\rho)\,\,\,\,\,\,\,\,\,\,\,\,\,\mathrm{and}\,\,\,\,\,\,\,\,\,\,\,\,\,Q_{\infty}\,=\,Q(\sigma\theta,\sigma\rho)\,.
\]
We also consider the family of boxes 
\[
\widetilde{Q}_{j}:=\,Q(\tilde{\theta}_{j},\tilde{\rho}_{j})\,,
\]
where, for $j\in\mathbb{N}_{0}$, 
\[
\tilde{\rho}_{j}:=\,\frac{\rho_{j}+\rho_{j+1}}{2}\,=\,\sigma\rho\,+\,\frac{3(1-\sigma)}{2^{j+2}}\,\rho\,,\,\,\,\,\,\,\,\,\,\,\,\,\,\tilde{\theta}_{j}:=\,\frac{\theta_{j}+\theta_{j+1}}{2}\,=\,\sigma\theta\,+\,\frac{3(1-\sigma)}{2^{j+2}}\,\theta\,.
\]
For these boxes, we have the inclusions
\begin{equation}
Q_{j+1}\subset\widetilde{Q}_{j}\subset Q_{j}\,,\,\,\,\,\,\,\,\,\,\,j\in\mathbb{N}_{0}\,.\label{eq:inclusions}
\end{equation}
We now introduce the sequence of increasing levels 
\begin{equation}
k_{j}:=\,k-\,\frac{k}{2^{j}}\,,\label{eq:levels}
\end{equation}
where $k$ is a positive number to be chosen later. We shall work
with inequality \eqref{eq:energy_est} written for the functions $(u-k_{j+1})_{+}$,
over the cylinders $Q_{j}$. The piecewise smooth cut-off function
$\zeta$ is taken to satisfy
\begin{equation}
\begin{cases}
\begin{array}{c}
{\displaystyle 0\leq\zeta\leq1\,,\,\,\,\,\,\,\,\,\,\,\zeta\equiv0\,\,\,\,\,\,\mathrm{on}\,\,\,\partial_{\mathrm{par}}Q_{j}\,,\,\,\,\,\,\,\,\,\,\,\zeta\equiv1\,\,\,\,\,\,\mathrm{in}\,\,\,\widetilde{Q}_{j}\,,}\vspace{4mm}\\
{\displaystyle \vert D\zeta\vert\leq\,\frac{2^{j+2}\,c}{(1-\sigma)\rho}\,\,,\,\,\,\,\,\,\,\,\,\,0\leq\,\partial_{t}\zeta\leq\,\frac{2^{j+2}\,c}{(1-\sigma)\theta}\,\,.\,\,\,\,\,\,\,\,\,\,\,\,\,\,\,\,\,\,\,\,\,\,}
\end{array}\end{cases}\label{eq:zeta}
\end{equation}
 With these choices, estimate \eqref{eq:energy_est} yields\begin{align}\label{eq:starting_point}
&\underset{-\theta_j\,<\,\tau\,<\,0}{\sup}\,\int_{K_{\rho_j}}(u-k_{j+1})_{+}^{2}\,\zeta^{p}(x,\tau)\,dx\,+\sum_{i=1}^{n}\iint_{Q_{j}\,\cap\,\{u\,>\,k_{j+1}\}}(\vert u_{x_{i}}\vert-\delta_{i})_{+}^{p}\,\zeta^{p}\,dx\,dt\nonumber\\
&\,\,\,\,\,\,\,\leq\,\frac{C_1\,2^{j}}{(1-\sigma)\theta}\iint_{Q_{j}}(u-k_{j+1})_{+}^{2}\,dx\,dt\,+\,\frac{C_1\,2^{jp}}{(1-\sigma)^{p}\,\rho^{p}}\iint_{Q_{j}}(u-k_{j+1})_{+}^{p}\,dx\,dt\,,
\end{align}where $C_{1}\equiv C_{1}(n,p)>0$. From definition (\ref{eq:levels}),
we immediately have 
\begin{equation}
\iint_{Q_{j}}(u-k_{j+1})_{+}^{p}\,dx\,dt\,\leq\iint_{Q_{j}}(u-k_{j})_{+}^{p}\,dx\,dt\,.\label{eq:banale}
\end{equation}
Now observe that, for all $s>0$, \begin{align}\label{eq:superlevel}
\iint_{Q_{j}}(u-k_{j})_{+}^{s}\,dx\,dt\,&\geq\iint_{Q_{j}\,\cap\,\{u\,>\,k_{j+1}\}}(u-k_{j})_{+}^{s}\,dx\,dt\nonumber\\
&\geq\,(k_{j+1}-k_{j})^{s}\,\vert A_{j+1}\vert\nonumber\\
&=\,\frac{k^{s}}{2^{(j+1)s}}\,\vert A_{j+1}\vert\,,
\end{align}where we have set 
\begin{equation}
\vert A_{j+1}\vert:=\,\mathrm{meas}\left\{ (x,t)\in Q_{j}:u(x,t)>k_{j+1}\right\} .\label{eq:measure}
\end{equation}
Then, using Hölder's inequality, (\ref{eq:banale}) and \eqref{eq:superlevel},
we get\begin{align}\label{eq:banale02}
\iint_{Q_{j}}(u-k_{j+1})_{+}^{2}\,dx\,dt\,&\leq\left(\iint_{Q_{j}}(u-k_{j+1})_{+}^{p}\,dx\,dt\right)^{\frac{2}{p}}\vert A_{j+1}\vert^{1-\frac{2}{p}}\nonumber\\
&\leq\,\frac{2^{(p-2)(j+1)}}{k^{p-2}}\iint_{Q_{j}}(u-k_{j})_{+}^{p}\,dx\,dt\,.
\end{align}Combining estimates \eqref{eq:starting_point}, (\ref{eq:banale})
and \eqref{eq:banale02}, and applying the properties $(\ref{eq:zeta})_{1}$
of $\zeta$, we obtain the following basic iterative inequalities:\begin{align}\label{eq:iterative00}
&\underset{-\tilde{\theta}_{j}\,<\,\tau\,<\,0}{\sup}\,\int_{K_{\tilde{\rho}_{j}}}(u(x,\tau)-k_{j+1})_{+}^{2}\,dx\,+\sum_{i=1}^{n}\iint_{\widetilde{Q}_{j}\,\cap\,\{u\,>\,k_{j+1}\}}(\vert u_{x_{i}}\vert-\delta_{i})_{+}^{p}\,dx\,dt\nonumber\\
&\,\,\,\,\,\,\,\leq\,\frac{C_1\,2^{jp}}{(1-\sigma)^{p}}\left(\frac{1}{\theta\,k^{p-2}}\,+\,\frac{1}{\rho^{p}}\right)\iint_{Q_{j}}(u-k_{j})_{+}^{p}\,dx\,dt\,.
\end{align}$\hspace*{1em}$To move forward, we construct a piecewise smooth cut-off
function $\tilde{\zeta}_{j}$ in $\widetilde{Q}_{j}$ such that
\[
\begin{cases}
\begin{array}{c}
{\displaystyle 0\leq\tilde{\zeta}_{j}\leq1\,,\,\,\,\,\,\,\,\,\,\,\tilde{\zeta}_{j}\equiv0\,\,\,\,\,\,\mathrm{on\,\,the\,\,lateral\,\,boundary\,\,of\,\,}\widetilde{Q}_{j},}\vspace{4mm}\\
\tilde{\zeta}_{j}\equiv1\,\,\,\,\,\,\mathrm{in}\,\,\,Q_{j+1}\,,\,\,\,\,\,\,\,\,\,\,{\displaystyle \vert D\tilde{\zeta}_{j}\vert\leq\,\frac{2^{j+2}}{(1-\sigma)\rho}\,.\,\,\,\,\,\,\,\,\,\,\,\,\,\,\,\,\,\,\,\,\,\,\,\,\,\,\,\,\,\,\,\,\,\,\,\,\,\,\,}
\end{array}\end{cases}
\]
Then the function $(u-k_{j+1})_{+}\,\tilde{\zeta}_{j}$ vanishes on
the lateral boundary of $\widetilde{Q}_{j}$ and, by Lemma \ref{lem:interpolation},
we have\begin{align}\label{eq:iterative1}
\iint_{\widetilde{Q}_{j}}(u-k_{j+1})_{+}^{q}\,\tilde{\zeta}_{j}^{q}\,dx\,dt\,\leq &\,\,C_{2}\left(\iint_{\widetilde{Q}_{j}}\vert D\,(u-k_{j+1})_{+}\vert^{p}\,dx\,dt\,+\iint_{\widetilde{Q}_{j}}(u-k_{j+1})_{+}^{p}\,\vert D\tilde{\zeta}_{j}\vert^{p}\,dx\,dt\right)\nonumber\\
&\,\times\left(\underset{-\tilde{\theta}_{j}\,<\,\tau\,<\,0}{\sup}\,\int_{K_{\tilde{\rho}_{j}}}(u(x,\tau)-k_{j+1})_{+}^{2}\,dx\right)^{\frac{p}{n}},
\end{align}where 
\begin{equation}
q:=\,{\displaystyle p\,\frac{n+2}{n}}\label{eq:q}
\end{equation}
and $C_{2}$ is a positive constant depending only on $n$ and $p$.\\
$\hspace*{1em}$At this point, we introduce the sequence of dimensionless
quantities
\begin{equation}
Y_{j}:=\,\tiltfiint_{Q_{j}}(u-k_{j})_{+}^{p}\,dx\,dt\,,\,\,\,\,\,\,\,\,\,\,j\in\mathbb{N}_{0}\,.\label{eq:Y_j}
\end{equation}
We shall derive an iterative inequality for $Y_{j}$ by estimating
the right-hand side of \eqref{eq:iterative1} by \eqref{eq:iterative00}.
Prior to this, by lengthy but elementary computations, we see that
\begin{equation}
\frac{\vert\widetilde{Q}_{j}\vert}{\vert Q_{j+1}\vert}\,<\,\left(\frac{3}{2}\right)^{n+1}\,\,\,\,\,\,\,\,\,\,\,\,\,\,\,\mathrm{and}\,\,\,\,\,\,\,\,\,\,\,\,\,\,\,\frac{\vert Q_{j}\vert}{\vert\widetilde{Q}_{j}\vert}\,<\,4^{n+1}\,.\label{eq:lenghty}
\end{equation}
Therefore, using the properties of $\tilde{\zeta}_{j}$, (\ref{eq:inclusions}),
(\ref{eq:lenghty}), Hölder's inequality, \eqref{eq:superlevel} with
$s=p$ and (\ref{eq:Y_j}), we obtain\begin{align*}
Y_{j+1}\,&=\,\tiltfiint_{Q_{j+1}}(u-k_{j+1})_{+}^{p}\,\tilde{\zeta}_{j}^{p}\,dx\,dt\,\leq\,\frac{\vert\widetilde{Q}_{j}\vert}{\vert Q_{j+1}\vert}\,\tiltfiint_{\widetilde{Q}_{j}}(u-k_{j+1})_{+}^{p}\,\tilde{\zeta}_{j}^{p}\,dx\,dt\\
&\leq\left(\frac{3}{2}\right)^{n+1}\vert\widetilde{Q}_{j}\vert^{\frac{p}{q}\,-\,1}\left(\tiltfiint_{\widetilde{Q}_{j}}(u-k_{j+1})_{+}^{q}\,\tilde{\zeta}_{j}^{q}\,dx\,dt\right)^{\frac{p}{q}}\vert A_{j+1}\vert^{1\,-\,\frac{p}{q}}\\
&\leq\,6^{n+1}\left(\tiltfiint_{\widetilde{Q}_{j}}(u-k_{j+1})_{+}^{q}\,\tilde{\zeta}_{j}^{q}\,dx\,dt\right)^{\frac{p}{q}}\left(\frac{\vert A_{j+1}\vert}{\vert Q_{j}\vert}\right)^{1\,-\,\frac{p}{q}}\\
&\leq\,C_{3}\left(\tiltfiint_{\widetilde{Q}_{j}}(u-k_{j+1})_{+}^{q}\,\tilde{\zeta}_{j}^{q}\,dx\,dt\right)^{\frac{p}{q}}\left(\frac{2^{jp}}{k^{p}}\,Y_{j}\right)^{1\,-\,\frac{p}{q}},
\end{align*}where $C_{3}\equiv C_{3}(n,p)>0$. We now estimate the last integral
via \eqref{eq:iterative1}, and subsequently estimate the right-hand
side of \eqref{eq:iterative1} using the inequality \eqref{eq:iterative00}.
Thus we obtain\begin{align}\label{eq:iterative_bis}
Y_{j+1}\,&\leq \,\frac{C_{3}}{\vert\widetilde{Q}_{j}\vert^{\frac{p}{q}}}\left(\iint_{\widetilde{Q}_{j}}\vert D\,(u-k_{j+1})_{+}\vert^{p}\,dx\,dt\,+\iint_{\widetilde{Q}_{j}}(u-k_{j+1})_{+}^{p}\,\vert D\tilde{\zeta}_{j}\vert^{p}\,dx\,dt\right)^{\frac{p}{q}}\nonumber\\
&\,\,\,\,\,\,\,\,\times\left(\underset{-\tilde{\theta}_{j}\,<\,\tau\,<\,0}{\sup}\,\int_{K_{\tilde{\rho}_{j}}}(u(x,\tau)-k_{j+1})_{+}^{2}\,dx\right)^{\frac{p^{2}}{nq}}\left(\frac{2^{jp}}{k^{p}}\,Y_{j}\right)^{1\,-\,\frac{p}{q}}\nonumber\\
&\leq\,\frac{C_{3}\,2^{\frac{p}{q}\left(\frac{p^{2}}{n}\,+\,q\,-\,p\right)j}}{(1-\sigma)^{\frac{p^{3}}{nq}}\,k^{\frac{p}{q}\,(q-p)}}\left(\iint_{\widetilde{Q}_{j}}\vert D\,(u-k_{j+1})_{+}\vert^{p}\,dx\,dt\,+\,\frac{2^{(j+2)p}}{(1-\sigma)^{p}\rho^{p}}\iint_{Q_{j}}(u-k_{j})_{+}^{p}\,dx\,dt\right)^{\frac{p}{q}}\nonumber\\
&\,\,\,\,\,\,\,\,\,\times\,\frac{\vert Q_{j}\vert^{\frac{p^{2}}{nq}}}{\vert\widetilde{Q}_{j}\vert^{\frac{p}{q}}}\left(\frac{1}{\theta\,k^{p-2}}\,+\,\frac{1}{\rho^{p}}\right)^{\frac{p^{2}}{nq}}Y_{j}^{\frac{p^{2}}{nq}\,+\,1\,-\,\frac{p}{q}}\nonumber\\
&\leq\,\frac{C_{3}\,2^{\frac{p}{q}\left(\frac{p^{2}}{n}\,+\,q\right)j}}{(1-\sigma)^{\frac{p^{3}}{nq}}\,k^{\frac{p}{q}\,(q-p)}}\,\,\frac{\vert Q_{j}\vert^{\frac{p^{2}}{nq}}}{\vert\widetilde{Q}_{j}\vert^{\frac{p}{q}}}\left(\frac{1}{\theta\,k^{p-2}}\,+\,\frac{1}{\rho^{p}}\right)^{\frac{p^{2}}{nq}}Z_{j}\,Y_{j}^{\frac{p^{2}}{nq}\,+\,1\,-\,\frac{p}{q}}\,,
\end{align}where, in the last line, we have set 
\begin{equation}
Z_{j}:=\left(\iint_{\widetilde{Q}_{j}}\vert D\,(u-k_{j+1})_{+}\vert^{p}\,dx\,dt\,+\,\frac{\vert Q_{j}\vert}{(1-\sigma)^{p}\rho^{p}}\,Y_{j}\right)^{\frac{p}{q}}.\label{eq:Z_j}
\end{equation}
Without loss of generality, we can now assume that $k\geq\rho$. Setting
\[
\delta:=\,\max\,\{\delta_{i}:i=1,\ldots,n\}\,,
\]
and using (\ref{eq:inclusions}), (\ref{eq:measure}), \eqref{eq:superlevel}
with $s=p$, \eqref{eq:iterative00}, (\ref{eq:Y_j}) and the fact
that $\frac{1}{k}\leq\frac{1}{\rho}$, we get\begin{align}\label{eq:iterative_ter}
Z_{j}\,&\leq\,C_{4}\left(\iint_{\widetilde{Q}_{j}}\sum_{i=1}^{n}\vert[(u-k_{j+1})_{+}]_{x_{i}}\vert^{p}\,dx\,dt\,+\,\frac{\vert Q_{j}\vert}{(1-\sigma)^{p}\rho^{p}}\,Y_{j}\right)^{\frac{p}{q}}\nonumber\\
&=\,C_{4}\left(\sum_{i=1}^{n}\iint_{\widetilde{Q}_{j}\,\cap\,\{u\,>\,k_{j+1}\}}\vert u_{x_{i}}\vert^{p}\,dx\,dt\,+\,\frac{\vert Q_{j}\vert}{(1-\sigma)^{p}\rho^{p}}\,Y_{j}\right)^{\frac{p}{q}}\nonumber\\
&\leq\,C_{4}\left(\sum_{i=1}^{n}\iint_{\widetilde{Q}_{j}\,\cap\,\{u\,>\,k_{j+1}\}}[(\vert u_{x_{i}}\vert-\delta_{i})_{+}+\delta_{i}]^{p}\,dx\,dt\,+\,\frac{\vert Q_{j}\vert}{(1-\sigma)^{p}\rho^{p}}\,Y_{j}\right)^{\frac{p}{q}}\nonumber\\
&\leq\,2^{\frac{p^{2}-p}{q}}C_{4}\left(\sum_{i=1}^{n}\iint_{\widetilde{Q}_{j}\,\cap\,\{u\,>\,k_{j+1}\}}(\vert u_{x_{i}}\vert-\delta_{i})_{+}^{p}\,dx\,dt\,+\,n\,\delta^{p}\,\vert A_{j+1}\vert\,+\,\frac{\vert Q_{j}\vert}{(1-\sigma)^{p}\rho^{p}}\,Y_{j}\right)^{\frac{p}{q}}\nonumber\\
&\leq\,\frac{C_{5}}{(1-\sigma)^{\frac{p^{2}}{q}}}\left[2^{jp}\left(\frac{1}{\theta\,k^{p-2}}\,+\,\frac{1}{\rho^{p}}\right)\vert Q_{j}\vert\,Y_{j}\,+\,\vert A_{j+1}\vert\,+\,\frac{1}{\rho^{p}}\,\vert Q_{j}\vert\,Y_{j}\right]^{\frac{p}{q}}\nonumber\\
&\leq\,\frac{C_{5}\,2^{(j+1)\frac{p^{2}}{q}}}{(1-\sigma)^{\frac{p^{2}}{q}}}\left[\left(\frac{1}{\theta\,k^{p-2}}\,+\,\frac{2}{\rho^{p}}\,+\,\frac{1}{k^{p}}\right)\vert Q_{j}\vert\,Y_{j}\right]^{\frac{p}{q}}\nonumber\\
&\leq\,\frac{C_{5}\,2^{(j+1)\frac{p^{2}}{q}}\,3^{\frac{p}{q}}\,}{(1-\sigma)^{\frac{p^{2}}{q}}}\left[\left(\frac{1}{\theta\,k^{p-2}}\,+\,\frac{1}{\rho^{p}}\right)\vert Q_{j}\vert\,Y_{j}\right]^{\frac{p}{q}},
\end{align}where $C_{4}\equiv C_{4}(n,p)>1$ and $C_{5}\equiv C_{5}(n,p,\delta)>1$.
Joining estimates \eqref{eq:iterative_bis} and \eqref{eq:iterative_ter},
we deduce 
\[
Y_{j+1}\,\leq\,\frac{C_{6}\,2^{\frac{p}{q}\left(\frac{p^{2}}{n}\,+\,q\,+\,p\right)j}}{(1-\sigma)^{\frac{p^{2}}{nq}\,(p+n)}\,k^{\frac{p}{q}\,(q-p)}}\,\,\frac{\vert Q_{j}\vert^{\frac{p}{nq}\,(p+n)}}{\vert\widetilde{Q}_{j}\vert^{\frac{p}{q}}}\,\left(\frac{1}{\theta\,k^{p-2}}\,+\,\frac{1}{\rho^{p}}\right)^{\frac{p}{nq}\,(p+n)}Y_{j}^{1\,+\,\frac{p^{2}}{nq}}
\]
for some positive constant $C_{6}$ depending only on $n$, $p$ and
$\delta$. From (\ref{eq:lenghty}) and (\ref{eq:inclusions}) again,
we obtain 
\[
\frac{\vert Q_{j}\vert^{\frac{p}{nq}\,(p+n)}}{\vert\widetilde{Q}_{j}\vert^{\frac{p}{q}}}\,\leq\left(\frac{4^{n+1}}{\vert Q_{j}\vert}\right)^{\frac{p}{q}}\vert Q_{j}\vert^{\frac{p}{nq}\,(p+n)}\,\leq\,4^{(n+1)\frac{p}{q}}\,\vert Q_{0}\vert^{\frac{p^{2}}{nq}}\,\leq\,C_{7}\,(\rho^{n}\theta)^{\frac{p^{2}}{nq}}\,,
\]
where $C_{7}\equiv C_{7}(n,p)>0$. Furthermore, we have 
\[
\left(\frac{1}{\theta\,k^{p-2}}\,+\,\frac{1}{\rho^{p}}\right)^{\frac{p}{nq}\,(p+n)}\leq\,C_{7}\left[\left(\frac{1}{\theta}\right)^{\frac{p+n}{p}}k^{(2-p)\,\frac{p+n}{p}}\,+\left(\frac{1}{\rho}\right)^{p+n}\right]^{\frac{p^{2}}{nq}}
\]
for a possibly different constant $C_{7}$. Finally, combining the
three previous estimates and recalling that $q:=p\,\frac{n+2}{n}$,
we arrive at the recursive inequalities 
\begin{equation}
Y_{j+1}\,\leq\,\frac{\widetilde{C}\,b^{j}}{(1-\sigma)^{p\,\frac{n+p}{n+2}}\,k^{\frac{2p}{n+2}}}\,\,\mathcal{A}_{k}^{\frac{p}{n+2}}\,Y_{j}^{1\,+\,\frac{p}{n+2}}\,,\label{eq:iter-ineq}
\end{equation}
where $\widetilde{C}$ is a positive constant depending only on $n$,
$p$ and $\delta$,
\begin{equation}
b:=\,2^{p\,\frac{p\,+\,2n\,+\,2}{n+2}}\,,\label{eq:b}
\end{equation}
\begin{equation}
\mathcal{A}_{k}:=\left(\frac{\rho^{p}}{\theta}\right)^{\frac{n}{p}}k^{(2-p)\,\frac{n+p}{p}}\,+\,\frac{\theta}{\rho^{p}}\,.\label{eq:A_k}
\end{equation}

\section{Proofs of the main results\label{sec:main proofs}}

\noindent $\hspace*{1em}$We next turn to the proof of Theorem \ref{thm:main}.
The argument relies on the iterative inequalities (\ref{eq:iter-ineq}),
which, combined with Lemma \ref{lem:Giusti}, yield the desired local
$L^{\infty}$ bounds (\ref{eq:superquad}) and (\ref{eq:quad}).

\noindent \begin{proof}[\bfseries{Proof of Theorem~\ref{thm:main}}]Let
us first consider the case $p>2$. After a translation, we may assume
that $(x_{0},t_{0})=(0,0)$. Therefore, we can make use of the iterative
inequalities (\ref{eq:iter-ineq}), where $Y_{j}$, $b$ and $\mathcal{A}_{k}$
are defined in (\ref{eq:Y_j}), (\ref{eq:b}) and (\ref{eq:A_k}),
respectively. Recalling that we assumed $k\geq\rho$ in obtaining
(\ref{eq:iter-ineq}), we now take $k$ so large that, of the two
terms composing $\mathcal{A}_{k}$, the second dominates the first,
i.e., 
\begin{equation}
k\,\ge\,\max\left\{ \rho,\left(\frac{\rho^{p}}{\theta}\right)^{\frac{1}{p-2}}\right\} .\label{eq:choice1}
\end{equation}
With this choice of $k$, we have 
\[
\mathcal{A}_{k}\,\leq\,\frac{2\theta}{\rho^{p}}\,.
\]
It follows from Lemma \ref{lem:Giusti} that $Y_{j}\rightarrow0$
as $j\rightarrow\infty$, provided we choose $k$ from 
\[
Y_{0}:=\,\tiltfiint_{Q(\theta,\rho)}u_{+}^{p}\,dx\,dt\,=\,C\,\,\frac{\rho^{p}}{\theta}\,(1-\sigma)^{n+p}\,k^{2}\,,
\]
where $C$ is a positive constant depending only on $n$, $p$ and
$\max\,\{\delta_{1},\ldots,\delta_{n}\}$. For such a choice and (\ref{eq:choice1}),
we obtain
\begin{equation}
\underset{Q(\sigma\theta,\sigma\rho)}{\mathrm{ess}\,\sup}\,\,u\,\leq\,\max\left\{ \rho,\left(\frac{\rho^{p}}{\theta}\right)^{\frac{1}{p-2}},\frac{C}{(1-\sigma)^{\frac{n+p}{2}}}\,\,\sqrt{\frac{\theta}{\rho^{p}}}\left(\tiltfiint_{Q(\theta,\rho)}u_{+}^{p}\,dx\,dt\right)^{\frac{1}{2}}\right\} ,\label{eq:sup_u}
\end{equation}
for a possibly different constant $C$. Now observe that $-u$ is
also a local weak solution of (\ref{eq:equation}). Then, replacing
$u$ with $-u$ in (\ref{eq:sup_u}), we get 
\[
\underset{Q(\sigma\theta,\sigma\rho)}{\mathrm{ess}\,\inf}\,\,u\,\geq\,-\max\left\{ \rho,\left(\frac{\rho^{p}}{\theta}\right)^{\frac{1}{p-2}},\frac{C}{(1-\sigma)^{\frac{n+p}{2}}}\,\,\sqrt{\frac{\theta}{\rho^{p}}}\left(\tiltfiint_{Q(\theta,\rho)}(-u)_{+}^{p}\,dx\,dt\right)^{\frac{1}{2}}\right\} .
\]
Combining the two previous estimates, we deduce the local $L^{\infty}$
bound in (\ref{eq:superquad}).\\
$\hspace*{1em}$Finally, if $p=2$ we have 
\[
\mathcal{A}_{k}=\left(\frac{\rho^{2}}{\theta}\right)^{\frac{n}{2}}+\,\frac{\theta}{\rho^{2}}\,,
\]
and arguing as above we obtain estimate (\ref{eq:quad}). This concludes
the proof.\end{proof}

\noindent $\hspace*{1em}$For the sake of completeness, we now detail
the modifications to the previous arguments that yield the

\noindent \begin{proof}[\bfseries{Proof of Corollary~\ref{cor:corollario}}]Let
$u$ be a local weak solution of (\ref{eq:AmbPass}) and let $\zeta$
be a piecewise smooth cut-off function chosen as in Section \ref{sec:energy_estimate}.
Now we set 
\[
A(\eta):=\,D_{\xi}G(\eta)\,,\,\,\,\,\,\eta\in\mathbb{R}^{n},
\]
where $G$ is the second function defined in (\ref{eq:FG}). Note
that 
\[
\vert A(Du)\vert=(\vert Du\vert-\lambda)_{+}^{p-1}\,\leq\,\vert Du\vert^{p-1}
\]
and
\[
\langle A(Du),Du\rangle=(\vert Du\vert-\lambda)_{+}^{p-1}\,\vert Du\vert\,\ge\,(\vert Du\vert-\lambda)_{+}^{p}\,.
\]
Therefore, we can proceed as in the proof of Proposition \ref{prop:PropEnergy},
thus obtaining, for every $k>0$, the following local energy estimate:\begin{align*}
&\underset{t_{0}-\theta\,<\,\tau\,<\,t_{0}}{\sup}\,\int_{[x_{0}+K_{\rho}]}(u-k)_{+}^{2}\,\zeta^{p}(x,\tau)\,dx\,+\iint_{[(x_{0},t_{0})+Q(\theta,\rho)]}(\vert Du\vert-\lambda)_{+}^{p}\,\zeta^{p}\,\mathds{1}_{\{u\,>\,k\}}\,dx\,dt\nonumber\\
&\,\,\,\,\,\,\,\leq\,p\iint_{[(x_{0},t_{0})+Q(\theta,\rho)]}(u-k)_{+}^{2}\,\zeta^{p-1}\,\partial_{t}\zeta\,dx\,dt\,+\,C(p)\iint_{[(x_{0},t_{0})+Q(\theta,\rho)]}(u-k)_{+}^{p}\,\vert D\zeta\vert^{p}\,dx\,dt\,.
\end{align*}Starting from this estimate, assuming again that $(x_{0},t_{0})=(0,0)$,
and using the same notations and arguments as in Section \ref{sec:IterIneq},
we find the iterative inequality\begin{align}\label{eq:iter_4th}
&\underset{-\tilde{\theta}_{j}\,<\,\tau\,<\,0}{\sup}\,\int_{K_{\tilde{\rho}_{j}}}(u(x,\tau)-k_{j+1})_{+}^{2}\,dx\,+\iint_{\widetilde{Q}_{j}\,\cap\,\{u\,>\,k_{j+1}\}}(\vert Du\vert-\lambda)_{+}^{p}\,dx\,dt\nonumber\\
&\,\,\,\,\,\,\,\leq\,\frac{C\,2^{jp}}{(1-\sigma)^{p}}\left(\frac{1}{\theta\,k^{p-2}}\,+\,\frac{1}{\rho^{p}}\right)\iint_{Q_{j}}(u-k_{j})_{+}^{p}\,dx\,dt\,,
\end{align}for every $j\in\mathbb{N}_{0}$. Moreover, we again arrive at estimate
\eqref{eq:iterative_bis}, where $q$, $Y_{j}$ and $Z_{j}$ are defined
respectively in (\ref{eq:q}), (\ref{eq:Y_j}) and (\ref{eq:Z_j}).
Without loss of generality, we can now assume that $k\geq\rho$. Thus,
using (\ref{eq:inclusions}), (\ref{eq:measure}), \eqref{eq:superlevel}
with $s=p$, \eqref{eq:iter_4th}, (\ref{eq:Y_j}) and the fact that
$\frac{1}{k}\leq\frac{1}{\rho}$, we get\begin{align}\label{eq:Zj_corol}
Z_{j}\,&=\left(\iint_{\widetilde{Q}_{j}\,\cap\,\{u\,>\,k_{j+1}\}}\vert Du\vert^{p}\,dx\,dt\,+\,\frac{\vert Q_{j}\vert}{(1-\sigma)^{p}\rho^{p}}\,Y_{j}\right)^{\frac{p}{q}}\nonumber\\
&\leq\left(\iint_{\widetilde{Q}_{j}\,\cap\,\{u\,>\,k_{j+1}\}}[(\vert Du\vert-\lambda)_{+}+\lambda]^{p}\,dx\,dt\,+\,\frac{\vert Q_{j}\vert}{(1-\sigma)^{p}\rho^{p}}\,Y_{j}\right)^{\frac{p}{q}}\nonumber\\
&\leq\,2^{\frac{p^{2}-p}{q}}\left(\iint_{\widetilde{Q}_{j}\,\cap\,\{u\,>\,k_{j+1}\}}(\vert Du\vert-\lambda)_{+}^{p}\,dx\,dt\,+\,\lambda^{p}\,\vert A_{j+1}\vert\,+\,\frac{\vert Q_{j}\vert}{(1-\sigma)^{p}\rho^{p}}\,Y_{j}\right)^{\frac{p}{q}}\nonumber\\
&\leq\,\frac{C_{1}}{(1-\sigma)^{\frac{p^{2}}{q}}}\left[2^{jp}\left(\frac{1}{\theta\,k^{p-2}}\,+\,\frac{1}{\rho^{p}}\right)\vert Q_{j}\vert\,Y_{j}\,+\,\vert A_{j+1}\vert\,+\,\frac{1}{\rho^{p}}\,\vert Q_{j}\vert\,Y_{j}\right]^{\frac{p}{q}}\nonumber\\
&\leq\,\frac{C_{1}\,2^{(j+1)\frac{p^{2}}{q}}}{(1-\sigma)^{\frac{p^{2}}{q}}}\left[\left(\frac{1}{\theta\,k^{p-2}}\,+\,\frac{2}{\rho^{p}}\,+\,\frac{1}{k^{p}}\right)\vert Q_{j}\vert\,Y_{j}\right]^{\frac{p}{q}}\nonumber\\
&\leq\,\frac{C_{1}\,2^{(j+1)\frac{p^{2}}{q}}\,3^{\frac{p}{q}}\,}{(1-\sigma)^{\frac{p^{2}}{q}}}\left[\left(\frac{1}{\theta\,k^{p-2}}\,+\,\frac{1}{\rho^{p}}\right)\vert Q_{j}\vert\,Y_{j}\right]^{\frac{p}{q}},
\end{align}where $C_{1}\equiv C_{1}(n,p,\lambda)>1$. Joining estimates \eqref{eq:iterative_bis}
and \eqref{eq:Zj_corol}, and arguing as in the final part of Section
\ref{sec:IterIneq}, we obtain the recursive inequalities (\ref{eq:iter-ineq}),
where $\widetilde{C}$ is now a positive constant depending only on
$n$, $p$ and $\lambda$, while $b$ and $\mathcal{A}_{k}$ are defined
in (\ref{eq:b}) and (\ref{eq:A_k}), respectively. The desired conclusion
then follows by proceeding exactly as in the proof of Theorem \ref{thm:main}.\end{proof}

\begin{singlespace}
\noindent \medskip{}

\noindent \textbf{Acknowledgments. }This work has been partially supported
by the INdAM\textminus GNAMPA 2025 Project ``Regolarità ed esistenza
per operatori anisotropi'' (CUP E5324001950001). The authors wish
to express their gratitude to the Department of Mathematics of the
University of Bologna. In addition, P. Ambrosio acknowledges financial
support under the National Recovery and Resilience Plan (NRRP), Mission
4, Component 2, Investment 1.1, Call for tender No. 104 published
on 2.2.2022 by the Italian Ministry of University and Research (MUR),
funded by the European Union - NextGenerationEU - Project PRIN\_CITTI
2022 - Title ``Regularity problems in sub-Riemannian structures''
- CUP J53D23003760006 - Bando 2022 - Prot. 2022F4F2LH.\bigskip{}
\textbf{Declarations.} On behalf of all authors, the corresponding
author states that there is no conflict of interest.
\end{singlespace}

\begin{singlespace}

\lyxaddress{\noindent \textbf{$\quad$}\\
$\hspace*{1em}$\textbf{Pasquale Ambrosio} \\
Dipartimento di Matematica, Università di Bologna\\
Piazza di Porta S. Donato 5, 40126 Bologna, Italy.\\
\textit{E-mail address}: pasquale.ambrosio@unibo.it}

\lyxaddress{\noindent $\hspace*{1em}$\textbf{Simone Ciani}\\
Dipartimento di Matematica, Università di Bologna\\
Piazza di Porta S. Donato 5, 40126 Bologna, Italy.\\
\textit{E-mail address}: simone.ciani3@unibo.it}
\end{singlespace}

\end{document}